# Gauge-measurable functions

Augusto C. PONCE and Jean VAN SCHAFTINGEN

Université catholique de Louvain, Institut de Recherche en Mathématique et Physique,
Chemin du Cyclotron 2, bte L7.01.02, 1348 Louvain-la-Neuve, Belgium

## Abstract

In 1973, E. J. McShane proposed an alternative definition of the Lebesgue integral based on Riemann sums, where gauges are used decide what tagged partitions are allowed. Such an approach does not require any preliminary knowledge of Measure theory. We investigate in this paper a definition of measurable functions also based on gauges. Its relation to the gauge-integrable functions that satisfy McShane's definition is obtained using elementary tools from Real Analysis. We show in particular a dominated integration property of gauge-measurable functions.

## 1 Introduction

In its classical original setting, the Lebesgue integral of a function $f$ is defined in terms of the outer Lebesgue measure of the measurable sublevel sets of $f$ [2, 11, 15]. Compared to the definition of the integral by Cauchy or Riemann as a limit of Riemann sums [4, 25], the definition of the Lebesgue integral seems somehow indirect: it is a limit of a sum of measures, where these measures are themselves computed as the infima or suprema of volumes.

This issue has led to the definition of *gauge integrals* as a way of recovering the original approach based on Riemann sums, without the defects associated to the Riemann integral of Riemann-integrable functions [20]. Around 1960, Kurzweil and Henstock independently defined a gauge integral which allows one to integrate more functions than the Lebesgue-integrable ones [10, 14]. A few years later, in 1973, McShane presented the Lebesgue integral itself as a gauge integral [22, 23]. We can rephrase McShane's definition as follows:

**Definition 1.1** (Gauge integrability). *A function $f : \mathbb{R}^d \to \mathbb{R}^p$ is* gauge-integrable *whenever there exists $I \in \mathbb{R}$ verifying the following property: for every $\epsilon > 0$, there exists a* *gauge $\gamma$ on $\mathbb{R}^d$ and a compact set $K \subset \mathbb{R}^d$ such that, for every finite set of disjoint rectangles $(R_i)_{i \in \{1,\ldots,k\}}$ in $\mathbb{R}^d$ that covers $K$ and every finite set of points $(c_i)_{i \in \{1,\ldots,k\}}$ in $\mathbb{R}^d$ satisfying*

$$R_i \subset \gamma(c_i) \quad \text{for every } i \in \{1, \ldots, k\},$$

*one has*

$$\Big| \sum_{i=1}^{k} f(c_i) \operatorname{vol}(R_i) - I \Big| \leq \epsilon.$$

In this definition, a *gauge* $\gamma$ on $\mathbb{R}^d$ is a function mapping each point of $x \in \mathbb{R}^d$ to an open set $\gamma(x) \subset \mathbb{R}^d$ such that $x \in \gamma(x)$; for example $\gamma(x)$ might be taken to be a non-empty open ball centered at $x$. A *rectangle* $R \subset \mathbb{R}^d$ is a set that can be written as $R = [a_1, b_1) \times \cdots \times [a_d, b_d)$, where $a_1 < b_1, \ldots, a_d < b_d$ are all real numbers; its volume is the positive number $\operatorname{vol}(R) = (b_1 - a_1) \cdots (b_d - a_d)$. Rectangles are *disjoint* whenever there intersection is empty, and the family $(R_i)_{i \in \{1,\ldots,k\}}$ *covers* $K$ if

$$\bigcup_{i=1}^{k} R_i \supset K.$$

The compact set $K$ corresponds in McShane's original definition to the complement of his gauge at infinity; the equivalent formulation above avoids compactifying the Euclidean space $\mathbb{R}^d$ and considering unbounded rectangles.

By Cousin's lemma, which is a variant of the Heine–Borel theorem, for any gauge $\gamma$ on $\mathbb{R}^d$ and any compact set $K \subset \mathbb{R}^d$, there always exist some finite set of disjoint rectangles $(R_i)_{i \in \{1,\ldots,k\}}$ that covers $K$ and points $(c_i)_{i \in \{1,\ldots,k\}}$ such that $R_i \subset \gamma(c_i)$ for every $i$ [23, Theorem IV-3-1]. This fact ensures the uniqueness of the integral $I$ of $f$, which entitles one to adopt the usual notation

$$\int_{\mathbb{R}^d} f := I.$$

A non-intuitive feature of the definition of the gauge integral above is that each tag $c_i$ need not belong to the rectangle $R_i$. Adding this restriction gives the broader definition of integral of Kurzweil and Henstock, which is a gauge definition of the Denjoy–Perron integral for which all derivatives of one-dimensional functions are integrable on bounded intervals [8, 17, 24]. This Kurzweil–Henstock integral has been taught by Jean Mawhin at the Université catholique de Louvain (UCL) for thirty years [18, 19], continuing the Louvain tradition of cutting-edge lectures on integration theory initiated by Ch.-J. de la Vallée Poussin with the Lebesgue integral at the beginning of the 20th century [5–7, 21]. The further restriction that the gauge $\gamma(x)$ contain some uniform ball $B_\delta(x)$ for some radius $\delta > 0$ independent of $x \in \mathbb{R}^d$ yields the classical Riemann integral.

Measurability of functions is not a prerequisite of McShane's definition of gauge integrability. This is an important aspect one should not neglect about the gauge integral that makes the Lebesgue integral readily available, without the need of any preliminary development of tools from Measure theory. This is an approach we have been pursuing at UCL since 2009.

When measurability is needed to state some integrability conditions, measurable functions have been defined as pointwise limits of integrable functions [23, Definition III-10-1] or almost everywhere limits of locally integrable step functions [1, §19; 16, Definition 3.5.3], or in terms of measurable sets whose characteristic function are, by definition, locally integrable functions [18, §6.B; 19, §13.7]. It thus seems that the straightforwardness of McShane's definition of the integral is lost in an *ad hoc* indirect definition of measurability based on the integral itself.

In order to remedy to this issue, we introduce here a direct definition of measurability of functions in terms of gauges inspired by Lusin's property for Lebesgue-measurable functions.

**Definition 1.2** (Gauge measurability). *A function $f : \mathbb{R}^d \to \mathbb{R}^p$ is gauge-measurable whenever, for every $\epsilon > 0$ and every $\eta > 0$, there exists a gauge $\gamma$ on $\mathbb{R}^d$ such that, for every finite set of disjoint rectangles $(R_i)_{i \in \{1,\dots,k\}}$ in $\mathbb{R}^d$ and finite sets of points $(c_i)_{i \in \{1,\dots,k\}}$ and $(c'_i)_{i \in \{1,\dots,k\}}$ in $\mathbb{R}^d$ satisfying*

$$|f(c_i) - f(c'_i)| \geq \eta \quad \text{and} \quad R_i \subset \gamma(c_i) \cap \gamma(c'_i) \quad \text{for every } i \in \{1,\dots,k\},$$

*one has*

$$\sum_{i=1}^{k} \text{vol}(R_i) \leq \epsilon.$$

The goal of this paper is to provide various properties of gauge-measurable functions that can be deduced using elementary ideas of Real Analysis. These are well-known properties of Lebesgue-measurable functions, and both notions of measurability are equivalent, but the main message we want to emphasize is that one can obtain these properties in a self-contained approach based on gauge integrability and gauge measurability. As an example, we show in Section 5 below that these two concepts are related through the following dominated integrability criterion for a function to be gauge-integrable:

**Theorem 1.3.** *A function $f : \mathbb{R}^d \to \mathbb{R}^p$ is gauge-integrable if and only if $f$ is gauge-measurable and there exists a gauge-integrable function $h : \mathbb{R}^d \to \mathbb{R}$ such that $|f| \leq h$ in $\mathbb{R}^d$.*

The paper is organized as follows. In Sections 2 and 3, we prove properties of gauge-measurable functions that can be straightforwardly obtained from the definition. Some of them will be superseded in later sections using two important properties of the gauge integral: the Absolute Cauchy criterion and the Dominated convergence theorem. In Section 4, we prove Lusin's theorem for gauge-measurable functions using an alternative formulation of the gauge measurability based on the inner measure of open sets in $\mathbb{R}^d$. We then prove Theorem 1.3 in Section 5. In Section 6, we prove the stability of gauge measurability under pointwise convergence. In Sections 7 and 8, we define gauge-measurable sets in the same spirit as for functions, and then we prove that every gauge-measurable function is the pointwise limit of gauge-measurable step functions. We thus recover the approach which leads to the Lebesgue integral.

## 2 Elementary properties

The goal of this section is to present some properties of gauge measurability that readily follows from its definition. We begin by noting that every continuous function is gauge-measurable.

**Proposition 2.1** (Gauge measurability of continuous functions)**.** *If the function $f : \mathbb{R}^d \to \mathbb{R}^p$ is continuous, then $f$ is gauge-measurable.*

*Proof.* Given a pair of points $c_i, c_i' \in \mathbb{R}^d$, by the triangle inequality for every $z \in \mathbb{R}^d$ we have
$$|f(c_i) - f(c_i')| \le |f(z) - f(c_i)| + |f(z) - f(c_i')|. \tag{2.1}$$
Using the continuity of $f$, we choose a gauge $\gamma$ in such a way that the right-hand side is always less than $\eta > 0$ provided that $\gamma(c_i) \cap \gamma(c_i') \ne \emptyset$. Indeed, given $\eta > 0$, for every $x \in \mathbb{R}^d$ we define
$$\gamma(x) = \Big\{ z \in \mathbb{R}^d \;\Big|\; |f(z) - f(x)| < \frac{\eta}{2} \Big\}.$$
In particular, $x \in \gamma(x)$; since the function $f$ is continuous, the set $\gamma(x)$ is open. If there exists $z \in \gamma(c_i) \cap \gamma(c_i')$, then by the choice of $\gamma$ we have simultaneously
$$|f(z) - f(c_i))| < \frac{\eta}{2} \quad \text{and} \quad |f(z) - f(c_i')| < \frac{\eta}{2}.$$
In view of (2.1), we then have
$$|f(c_i) - f(c_i')| < \eta.$$

Therefore, no matter what $\epsilon > 0$ we take, there is no finite family of rectangles $(R_i)_{i \in \{1,\dots,k\}}$ that needs to be checked in Definition 1.2, so the latter is automatically satisfied by the continuous function $f$. □

**Proposition 2.2** (Composition with uniformly continuous functions)**.** *If the function $f : \mathbb{R}^d \to \mathbb{R}^p$ is gauge-measurable and the function $\Phi : \mathbb{R}^p \to \mathbb{R}^\ell$ is uniformly continuous, then the composition $\Phi \circ f : \mathbb{R}^d \to \mathbb{R}^\ell$ is gauge-measurable.*

This property is reminiscent of the integrability of compositions with Lipschitz functions for the gauge integral [23, Theorem II-3-1]; the class of admissible functions is larger here because the gauge measurability is a more qualitative property than the gauge integrability. As for the gauge integrability, Proposition 2.2 is not the end of the story: we will prove is Section 6 with more elaborate tools that the proposition remains true when the function $g$ is merely continuous; see Proposition 6.5 below.

*Proof of Proposition 2.2.* Given $\eta > 0$, by definition of uniform continuity there exists $\delta > 0$ such that, for every $y, z \in \mathbb{R}^d$ satisfying $|y - z| < \delta$, one has $|\Phi(y) - \Phi(z)| < \eta$. This is equivalent to saying that if $|\Phi(y) - \Phi(z)| \geq \eta$, then $|y - z| \geq \delta$. Hence, for every pair of points $c_i, c_i' \in \mathbb{R}^d$ such that

$$|(\Phi \circ f)(c_i) - (\Phi \circ f)(c_i')| \geq \eta, \tag{2.2}$$

we have

$$|f(c_i) - f(c_i')| \geq \delta. \tag{2.3}$$

Given $\epsilon > 0$, by Definition 1.2 of gauge measurability of $f$ with parameter $\delta$ there exists a gauge $\gamma$ on $\mathbb{R}^d$ such that, for every finite set of disjoint rectangles $(R_i)_{i \in \{1,\dots,k\}}$ and finite sets of points $(c_i)_{i \in \{1,\dots,k\}}$ and $(c_i')_{i \in \{1,\dots,k\}}$ in $\mathbb{R}^d$ satisfying (2.2) and $R_i \subset \gamma(c_i) \cap \gamma(c_i')$ for every $i$, we have that (2.3) also holds for every $i$, and then by the choice of the gauge $\gamma$,

$$\sum_{i=1}^{k} \mathrm{vol}\,(R_i) \leq \epsilon.$$

The function $\Phi \circ f$ is thus gauge-measurable. $\square$

An interesting consequence of Proposition 2.2 is that the family of gauge-measurable functions forms a vector space, and the product of two bounded gauge-measurable functions is also gauge-measurable. We provide an independent proof of these facts in the next section for the sake of clarity. The latter property concerning the product will be superseded later on by using the measurable stability under pointwise convergence, which allows one to remove the boundedness assumption of the functions; see Corollary 6.4. For the moment, we shall restrict ourselves to the case of uniform limits of gauge-measurable functions:

**Proposition 2.3** (Uniform limit). *Let $(f_n)_{n \in \mathbb{N}}$ be a sequence of gauge-measurable functions from $\mathbb{R}^d$ to $\mathbb{R}^p$. If the sequence $(f_n)_{n \in \mathbb{N}}$ converges uniformly to the function $f : \mathbb{R}^d \to \mathbb{R}^p$, then $f$ is gauge-measurable.*

*Proof.* For every pair of points $c_i, c_i' \in \mathbb{R}^d$ and $n \in \mathbb{N}$, by the triangle inequality we have

$$|f_n(c_i) - f_n(c_i')| \geq |f(c_i) - f(c_i')| - |f_n(c_i) - f(c_i)| - |f_n(c_i') - f(c_i')|.$$

Given $\eta > 0$, by the definition of uniform convergence there exists $n \in \mathbb{N}$ such that, for every $x \in \mathbb{R}^d$, $|f_n(x) - f(x)| \leq \eta/4$. Hence, assuming that

$$|f(c_i) - f(c_i')| \geq \eta, \tag{2.4}$$

we have

$$|f_n(c_i) - f_n(c_i')| \geq \frac{\eta}{2}. \tag{2.5}$$

Given $\epsilon > 0$, let $\gamma$ be a gauge on $\mathbb{R}^d$ given by the definition of gauge measurability of $f_n$ with parameter $\eta/2$. For every finite set of disjoint rectangles $(R_i)_{i \in \{1,\dots,k\}}$ and points $(c_i)_{i \in \{1,\ldots,k\}}$ and $(c'_i)_{i \in \{1,\ldots,k\}}$ satisfying (2.4) and $R_i \subset \gamma(c_i) \cap \gamma(c'_i)$ for every $i$, we then have that (2.5) is satisfied by $f_n$ for every $i$, and then, by the choice of $\gamma$,

$$\sum_{i=1}^{k} \mathrm{vol}\,(R_i) \leq \epsilon.$$

The function $f$ is thus gauge-integrable. □

## 3 Algebraic stability

We show that the class of gauge-measurable functions forms a vector space:

**Proposition 3.1** (Linearity)**.** *If the functions $f : \mathbb{R}^d \to \mathbb{R}^p$ and $g : \mathbb{R}^d \to \mathbb{R}^p$ are gauge-measurable and $\lambda \in \mathbb{R}$, then $f + g$ and $\lambda f$ are gauge-measurable.*

*Proof.* We focus on the proof that the function $f + g$ is gauge measurable; the case of $\lambda f$ is left as an exercise (see also Proposition 2.2). For every pair of points $c_i, c'_i \in \mathbb{R}^d$, by the triangle inequality we have

$$|(f+g)(c_i) - (f+g)(c'_i)| \leq |f(c_i) - f(c'_i)| + |g(c_i) - g(c'_i)|.$$

Given $\eta > 0$, and assuming that

$$|(f+g)(c_i) - (f+g)(c'_i)| \geq \eta, \tag{3.1}$$

then we necessarily have

$$|f(c_i) - f(c'_i)| \geq \frac{\eta}{2} \quad \text{or} \quad |g(c_i) - g(c'_i)| \geq \frac{\eta}{2}. \tag{3.2}$$

Given $\epsilon > 0$, let $\gamma_1$ and $\gamma_2$ be two gauges on $\mathbb{R}^d$ arising from the definitions of gauge measurability of $f$ and $g$, respectively, with parameters $\epsilon/2$ and $\eta/2$. Consider the gauge $\gamma$ defined for $x \in \mathbb{R}^d$ by $\gamma(x) = \gamma_1(x) \cap \gamma_2(x)$. For a finite collection of disjoint rectangles $(R_i)_{i \in \{1,\ldots,k\}}$ and finite sets of points $(c_i)_{i \in \{1,\ldots,k\}}$ and $(c'_i)_{i \in \{1,\ldots,k\}}$ in $\mathbb{R}^d$ verifying (3.1) and $R_i \subset \gamma(c_i) \cap \gamma(c'_i)$ for every $i \in \{1,\ldots,k\}$, let us denote by $I_1$ the set of indices $i$ for which the first inequality in (3.2) holds for $f$ and by $I_2$ the set of indices $i$ for which the second inequality in (3.2) holds for $g$. We can thus assert that

$$\{1,\ldots,k\} = I_1 \cup I_2. \tag{3.3}$$

We have in particular $R_i \subset \gamma_1(c_i) \cap \gamma_1(c'_i)$ for every $i \in I_1$, and thus by the choice of $\gamma_1$,

$$\sum_{i \in I_1} \mathrm{vol}\,(R_i) \leq \frac{\epsilon}{2}.$$

We also have $R_i \subset \gamma_2(c_i) \cap \gamma_2(c'_i)$ for every $i \in I_2$, and thus by the choice of $\gamma_2$,

$$\sum_{i \in I_2} \mathrm{vol}\,(R_i) \leq \frac{\epsilon}{2}.$$

Since the sets $I_1$ and $I_2$ cover $\{1, \ldots, k\}$, we deduce that

$$\sum_{i=1}^{k} \operatorname{vol}(R_i) \leq \frac{\epsilon}{2} + \frac{\epsilon}{2} = \epsilon.$$

Therefore, the function $f + g$ is gauge-measurable. □

Using a similar idea, one shows that the product of bounded gauge-measurable functions is also gauge-measurable. The conclusion is still true without assuming the functions are bounded, but the proof is more subtle and will be presented in Section 6 below.

**Proposition 3.2** (Product of bounded functions). *If the functions $f : \mathbb{R}^d \to \mathbb{R}^p$ and $g : \mathbb{R}^d \to \mathbb{R}$ are gauge-measurable and bounded, then $fg$ is also gauge-measurable.*

*Proof.* Take $M > 0$ and $N > 0$ such that $|f| \leq M$ and $|g| \leq N$ in $\mathbb{R}^d$. Given finite sets of points $(c_i)_{i \in \{1,\ldots,k\}}$ and $(c'_i)_{i \in \{1,\ldots,k\}}$ in $\mathbb{R}^d$, by the triangle inequality for every $x \in \mathbb{R}^d$ we have

$$\begin{aligned}|(fg)(c_i) - (fg)(c'_i)| &\leq |f(c_i) - f(c'_i)|\,|g(c_i)| + |f(c'_i)|\,|g(c_i) - g(c'_i)| \\ &\leq N|f(c_i) - f(c'_i)| + M|g(c_i) - g(c'_i)|.\end{aligned}$$

Given $\eta > 0$, if for every $i \in \{1, \ldots, k\}$ we have

$$|(fg)(c_i) - (fg)(c'_i)| \geq \eta,$$

then necessarily

$$|f(c_i) - f(c'_i)| \geq \frac{\eta}{2N} \quad \text{or} \quad |g(c_i) - g(c'_i)| \geq \frac{\eta}{2M}.$$

As in the previous proof, one defines the subsets of indices $I_1$ and $I_2$ accordingly, so that the counterpart of (3.3) also holds in this case. One can now proceed along the lines of the proof of Proposition 3.1 to deduce that $fg$ is gauge-measurable. □

## 4 Lusin's theorem

We now relate the notion of gauge measurability with Lusin's theorem, which trivially extends Proposition 2.1 for continuous functions:

**Proposition 4.1** (Lusin's theorem). *A function $f : \mathbb{R}^d \to \mathbb{R}^p$ is gauge-measurable if and only if, for every $\epsilon > 0$, there exists a closed set $C \subset \mathbb{R}^d$ such that the restriction $f|_C$ is continuous and the inner measure of the open set $\mathbb{R}^d \setminus C$ satisfies $\mu(\mathbb{R}^d \setminus C) \leq \epsilon$.*

We recall the notion of inner measure of an open set $U \subset \mathbb{R}^d$:

$$\mu(U) := \sup \Big\{\sum_{i=1}^{k} \text{vol}(R_i) \,\Big|\, (R_i)_{i \in \{1,\ldots,k\}} \text{ is a finite disjoint family of rectangles} \\ \text{such that } \bar{R}_i \subset U \text{ for every } i \in \{1,\ldots,k\}\Big\}.$$

Observe that $\mu$ is nondecreasing and countably subadditive.

Lusin's theorem above gives the equivalence between gauge measurability and the measurability in the sense of Bourbaki, defined in terms of Lusin's property [3, Definition IV-§5-1]. To prove Proposition 4.1 above, we rely on the following lemma which reformulates Definition 1.2 without relying on tagged partitions:

**Lemma 4.2** (Gauge-intersection characterization). *The function $f : \mathbb{R}^d \to \mathbb{R}^p$ is gauge-measurable if and only if, for every $\epsilon > 0$ and every $\eta > 0$, there exists a gauge $\gamma$ on $\mathbb{R}^d$ such that the open set*

$$U_{\gamma,\eta} := \bigcup_{\substack{x,z \in \mathbb{R}^d \\ |f(x)-f(z)| \geq \eta}} \big(\gamma(x) \cap \gamma(z)\big)$$

*satisfies $\mu(U_{\gamma,\eta}) \leq \epsilon$.*

A byproduct of Lemma 4.2 is the invariance of gauge measurability under bi-Lipschitz homeomorphisms of $\mathbb{R}^d$, which includes isometries.

*Proof of Lemma 4.2.* "$\Longleftarrow$". Given $\eta > 0$ and a gauge $\gamma$ on $\mathbb{R}^d$, take a finite disjoint family of rectangles $(R_i)_{i \in \{1,\ldots,k\}}$ and finite sets of points $(c_i)_{i \in \{1,\ldots,k\}}$ and $(c'_i)_{i \in \{1,\ldots,k\}}$ in $\mathbb{R}^d$ such that

$$|f(c_i) - f(c'_i)| \geq \eta \quad \text{and} \quad R_i \subset \gamma(c_i) \cap \gamma(c'_i) \quad \text{for every } i.$$

In particular, $R_i \subset \gamma(c_i) \cap \gamma(c'_i) \subset U_{\gamma,\eta}$, hence by definition of the inner measure $\mu(U_{\gamma,\eta})$ we have

$$\sum_{i=1}^{k} \text{vol}(R_i) \leq \mu(U_{\gamma,\eta}).$$

To conclude it suffices to choose the gauge $\gamma$ so that, for any given $\epsilon > 0$, we have $\mu(U_{\gamma,\eta}) \leq \epsilon$.

"$\Longrightarrow$". Assume that the function $f$ is gauge-measurable, and let $\gamma$ be a gauge on $\mathbb{R}^d$ given by Definition 1.2 for some $\epsilon > 0$ and $\eta > 0$. If $(R_i)_{i \in \{1,\ldots,k\}}$ is a finite disjoint family of rectangles such that $\bar{R}_i \subset U_{\gamma,\eta}$ for every $i$, then, by compactness of $\bar{R}_i$, the rectangle $R_i$ can be covered by a finite collection of sets of the form $\gamma(x) \cap \gamma(z)$ such that $x, z \in \mathbb{R}^d$ and $|f(x) - f(z)| \geq \eta$. By a suitable subdivision of the rectangles $(R_i)_{i \in \{1,\ldots,k\}}$ into smaller rectangles, which does not change their total volume, we can thus assume without loss of generality that, for every $i \in \{1,\ldots,k\}$, there exist points $x, z \in \mathbb{R}^d$ such that

$$R_i \subset \gamma(x) \cap \gamma(z) \quad \text{and} \quad |f(x) - f(z)| \geq \eta.$$

[Such a subdivision is allowed since the points $x$ and $z$ are not required to belong to $R_i$.] We then choose $c_i = x$ and $c'_i = z$. The finite sets of points $(c_i)_{i\in\{1,\ldots,k\}}$ and $(c'_i)_{i\in\{1,\ldots,k\}}$ satisfy the conditions of Definition 1.2, and we deduce that

$$\sum_{i=1}^{k} \mathrm{vol}\,(R_i) \leq \epsilon.$$

Since the family of rectangles $(R_i)_{i\in\{1,\ldots,k\}}$ is chosen arbitrarily, we thus have that $\mu(U_{\gamma,\eta}) \leq \epsilon$. □

*Proof of Proposition 4.1.* We first observe that given $\eta > 0$, a gauge $\gamma$ on $\mathbb{R}^d$, and $z \in \mathbb{R}^d$, then for every $x \in \gamma(z) \setminus U_{\gamma,\eta}$ we have

$$|f(x) - f(z)| < \eta. \tag{4.1}$$

Indeed, since $x \in \gamma(x) \cap \gamma(z)$ and $x \notin U_{\gamma,\eta}$, the set $\gamma(x) \cap \gamma(z)$ is not contained in $U_{\gamma,\eta}$, hence $x$ and $z$ are not admissible indices in the union that defines the set $U_{\gamma,\eta}$. We deduce that (4.1) holds.

Proceeding with the proof of the proposition, we now assume that the function $f$ is gauge-measurable and let $\epsilon > 0$. For each $n \in \mathbb{N}$, by Lemma 4.2 there exists a gauge $\gamma_n$ on $\mathbb{R}^d$ such that

$$\mu\bigl(U_{\gamma_n, 1/2^n}\bigr) \leq \frac{\epsilon}{2^{n+1}}.$$

We set $C = \mathbb{R}^d \setminus \bigcup_{n\in\mathbb{N}} U_{\gamma_n, 1/2^n}$. By countable subadditivity of $\mu$, we have

$$\mu(\mathbb{R}^d \setminus C) \leq \sum_{n\in\mathbb{N}} \mu\bigl(U_{\gamma_n, 1/2^n}\bigr) \leq \sum_{n\in\mathbb{N}} \frac{\epsilon}{2^{n+1}} = \epsilon.$$

It remains to prove that the restricted function $f|_C$ is continuous at any point $z \in C$. For every $x \in \gamma_n(z) \cap C \subset \gamma_n(z) \setminus U_{\gamma_n, 1/2^n}$, we deduce from estimate (4.1) above that

$$|f(x) - f(z)| < \frac{1}{2^n}.$$

Since this estimate holds on the relatively open subset $\gamma_n(z) \cap C$ of $C$ and $n \in \mathbb{N}$ is arbitrary, we deduce that the function $f|_C$ is continuous at $z$.

Conversely, take a closed set $C$ such that the restriction $f|_C$ is continuous. For every $\eta > 0$, the set

$$\gamma(x) = \mathbb{R}^d \setminus \left\{ w \in C \,\bigl|\, |f(x) - f(w)| \geq \frac{\eta}{2} \right\}$$

contains $x$ and is open in $\mathbb{R}^d$, since the function $f|_C$ is continuous and the set $C$ is closed. Hence, $\gamma$ is a gauge on $\mathbb{R}^d$. We now observe that if $x, z \in \mathbb{R}^d$ and $|f(x) - f(z)| \geq \eta$, then

$$\gamma(x) \cap \gamma(z) \cap C = \emptyset.$$

Indeed, if this were not true, there would exist $w \in \gamma(x) \cap \gamma(z) \cap C$. Since $w \in C$, we would have, by definition of $\gamma$, $|f(x) - f(w)| < \eta/2$ and $|f(z) - f(w)| < \eta/2$ and thus by the triangle inequality $|f(x) - f(z)| < \eta$, which would be a contradiction.

We thus have $U_{\gamma,\eta} \subset \mathbb{R}^d \setminus C$, and then by monotonicity of the inner measure $\mu$,

$$\mu(U_{\gamma,\eta}) \leq \mu(\mathbb{R}^d \setminus C).$$

Given $\epsilon > 0$, by the Lusin property satisfied by the function $f$, we may choose the closed set $C$ so as to have $\mu(\mathbb{R}^d \setminus C) \leq \epsilon$. We conclude from Lemma 4.2 that the function $f$ is gauge-measurable. □

## 5 Gauge measurability and integrability

The goal of this section is to establish Theorem 1.3. The relationship between gauge measurability and gauge integrability relies on the following Absolute Cauchy criterion for gauge-integrable functions [23, Theorem II-2-4] (see also [13, Lemma 5.13]).

**Proposition 5.1** (Absolute Cauchy criterion). *The function $f : \mathbb{R}^d \to \mathbb{R}^p$ is gauge-integrable if and only, for every $\epsilon > 0$, there exist a gauge $\gamma$ on $\mathbb{R}^d$ and a compact subset $K \subset \mathbb{R}^d$ such that the following properties hold:*

(i) *for every finite set of disjoint rectangles $(R_i)_{i \in \{1,...,k\}}$ in $\mathbb{R}^d$ and every finite sets of points $(c_i)_{i \in \{1,...,k\}}$ and $(c'_i)_{i \in \{1,...,k\}}$ satisfying $R_i \subset \gamma(c_i) \cap \gamma(c'_i)$ for every $i$, one has*

$$\sum_{i=1}^{k} |f(c_i) - f(c'_i)| \operatorname{vol}(R_i) \leq \epsilon.$$

(ii) *for every finite set of disjoint rectangles $(R_i)_{i \in \{1,...,k\}}$ in $\mathbb{R}^d \setminus K$ and every finite set of points $(c_i)_{i \in \{1,...,k\}}$ such that $R_i \subset \gamma(c_i)$ for every $i$, one has*

$$\sum_{i=1}^{k} |f(c_i)| \operatorname{vol}(R_i) \leq \epsilon.$$

This condition is a Cauchy criterion because it is an integrability criterion that does not require nor gives the value of the integral. It is an *absolute* Cauchy condition because the norm is taken inside the Riemann sum. An important consequence of Proposition 5.1 is the fact that if $f : \mathbb{R}^d \to \mathbb{R}^p$ is gauge-integrable and if $\Phi$ is a Lipschitz function such that $\Phi(0) = 0$, then the composite function $\Phi \circ f$ is also gauge-integrable [23, Theorem II-3-1]. In particular, $|f|$ is gauge-integrable whenever $f$ is gauge-integrable.

We first consider the question of the gauge measurability of gauge-integrable functions.

**Proposition 5.2** (Gauge measurability). *If $f : \mathbb{R}^d \to \mathbb{R}^p$ is gauge-integrable, then $f$ is gauge-measurable.*

*Proof.* Let $\eta > 0$ and take a finite set of disjoint rectangles $(R_i)_{i \in \{1,\dots,k\}}$ and finite sets of points $(c_i)_{i \in \{1,\dots,k\}}$ and $(c'_i)_{i \in \{1,\dots,k\}}$ in $\mathbb{R}^d$ such that

$$|f(c_i) - f(c'_i)| \geq \eta \quad \text{for every } i.$$

Then, we have

$$\sum_{i=1}^{k} \operatorname{vol}(R_i) \leq \frac{1}{\eta} \sum_{i=1}^{k} |f(c_i) - f(c'_i)| \operatorname{vol}(R_i). \tag{5.1}$$

Applying Property $(i)$ of the Absolute Cauchy criterion with parameter $\eta\epsilon$, there exists a gauge $\gamma$ on $\mathbb{R}^d$ such that if $R_i \subset \gamma(c_i) \cap \gamma(c'_i)$, then the sum in the right-hand side of (5.1) is smaller than $\eta\epsilon$, and we get

$$\sum_{i=1}^{k} \operatorname{vol}(R_i) \leq \frac{1}{\eta} \cdot \eta\epsilon = \epsilon.$$

We deduce that the function $f$ is gauge-measurable in view of Definition 1.2. $\square$

We now handle the reverse implication of Theorem 1.3 under the additional assumption that $f$ is a bounded function.

**Proposition 5.3** (Dominated integrability for bounded functions). *If $f : \mathbb{R}^d \to \mathbb{R}^p$ is gauge-measurable and bounded and if $|f| \leq h$ in $\mathbb{R}^d$ for some gauge-integrable function $h : \mathbb{R}^d \to \mathbb{R}$, then $f$ is gauge-integrable.*

*Proof.* Property $(ii)$ of the Absolute Cauchy criterion is satisfied by $h$, hence also by $f$. We now focus on Property $(i)$. For this purpose, let $(R_i)_{i \in \{1,\dots,k\}}$ be a finite collection of disjoint rectangles, and let $(c_i)_{i \in \{1,\dots,k\}}$ and $(c'_i)_{i \in \{1,\dots,k\}}$ be finitely many points in $\mathbb{R}^d$. Given $\eta > 0$ and a compact subset $K \subset \mathbb{R}^d$, we can relabel the rectangles and points simultaneously so as to have

(a) for every $i \in \{1,\dots,m\}$, $|f(c_i) - f(c'_i)| \geq \eta$,

(b) for every $i \in \{m+1,\dots,l\}$, $|f(c_i) - f(c'_i)| < \eta$ and $R_i \cap K \neq \emptyset$,

(c) for every $i \in \{l+1,\dots,k\}$, $|f(c_i) - f(c'_i)| < \eta$ and $R_i \cap K = \emptyset$,

for some integers $0 \leq m \leq l \leq k$; some of these condition might be empty, and in this case one simply ignores them.

By the assumption of boundedness of $f$, there exists $M > 0$ such that, for every $x \in \mathbb{R}^d$, $|f(x)| \leq M$. By the triangle inequality, we then have

$$\sum_{i=1}^{m} |f(c_i) - f(c'_i)| \operatorname{vol}(R_i) \leq 2M \sum_{i=1}^{m} \operatorname{vol}(R_i),$$

and, by $(b)$,
$$\sum_{i=m+1}^{l} |f(c_i) - f(c_i')| \operatorname{vol}(R_i) \leq \eta \sum_{i=m+1}^{l} \operatorname{vol}(R_i).$$

Since $|f| \leq h$ in $\mathbb{R}^d$, we also have
$$\sum_{i=l+1}^{k} |f(c_i) - f(c_i')| \operatorname{vol}(R_i) \leq \sum_{i=l+1}^{k} \big(|f(c_i)| + |f(c_i')|\big) \operatorname{vol}(R_i)$$
$$\leq \sum_{i=l+1}^{k} \big(h(c_i) + h(c_i')\big) \operatorname{vol}(R_i).$$

These are the three main estimates that we need in the sequel. We now proceed to choose the gauge $\gamma$ that yields the Absolute Cauchy criterion for $f$.

Given $\epsilon > 0$, by Property $(ii)$ of the Absolute Cauchy criterion satisfied by $h$ with parameter $\epsilon/6$, we can take the compact set $K \subset \mathbb{R}^d$ and a gauge $\gamma_1$ on $\mathbb{R}^d$ such that if $R_i \subset \gamma_1(c_i) \cap \gamma_1(c_i')$ for every $i \in \{l+1, \ldots, k\}$, then we have
$$\sum_{i=l+1}^{k} |f(c_i) - f(c_i')| \operatorname{vol}(R_i) \leq \sum_{i=l+1}^{k} \big(h(c_i) + h(c_i')\big) \operatorname{vol}(R_i) \leq \frac{\epsilon}{6} + \frac{\epsilon}{6} = \frac{\epsilon}{3}.$$

Fix a bounded open set $U \subset \mathbb{R}^d$ that contains $K$, and take the gauge $\gamma_2$ on $\mathbb{R}^d$ defined by $\gamma_2(x) = U$ if $x \in U$ and $\gamma_2(x) = \mathbb{R}^d \setminus K$ if $x \notin U$. Observe that if $R_i \subset \gamma_2(c_i) \cap \gamma_2(c_i')$ for every $i \in \{m+1, \ldots, l\}$, then since $R_i \cap K \neq \emptyset$, we necessarily have $\gamma_2(c_i) = \gamma_2(c_i') = U$, and thus $R_i \subset U$. By the definition of the inner measure $\mu$, and choosing $\eta > 0$ so as to have $\eta \mu(U) \leq \epsilon/3$, we then get
$$\sum_{i=m+1}^{l} |f(c_i) - f(c_i')| \operatorname{vol}(R_i) \leq \eta \sum_{i=m+1}^{l} \operatorname{vol}(R_i) \leq \eta \mu(U) \leq \frac{\epsilon}{3}.$$

By the definition of gauge measurability of $f$ with $\epsilon/6M$ and $\eta$ chosen as above, there exists a gauge $\gamma_3$ on $\mathbb{R}^d$ such that if $R_i \subset \gamma_3(c_i) \cap \gamma_3(c_i')$ for every $i \in \{1, \ldots, m\}$, then we have
$$\sum_{i=1}^{m} |f(c_i) - f(c_i')| \operatorname{vol}(R_i) \leq 2M \sum_{i=1}^{m} \operatorname{vol}(R_i) \leq 2M \cdot \frac{\epsilon}{6M} = \frac{\epsilon}{3}.$$

Combining these three estimates, we get
$$\sum_{i=1}^{k} |f(c_i) - f(c_i')| \operatorname{vol}(R_i) \leq \frac{\epsilon}{3} + \frac{\epsilon}{3} + \frac{\epsilon}{3} = \epsilon,$$

and thus $f$ satisfies the Absolute Cauchy criterion with the gauge $\gamma$ defined for $x \in \mathbb{R}^d$ by $\gamma(x) = \gamma_1(x) \cap \gamma_2(x) \cap \gamma_3(x)$. Hence, $f$ is gauge-integrable by Proposition 5.1. $\square$

The boundedness assumption of $f$ can be removed using the Dominated convergence theorem for gauge-integrable functions [23, Theorem II-10-1]:

**Proposition 5.4** (Dominated convergence). *Let $(f_n)_{n \in \mathbb{N}}$ be a sequence of gauge-integrable functions from $\mathbb{R}^d$ to $\mathbb{R}^p$. If $(f_n)_{n \in \mathbb{N}}$ converges pointwise to the function $f : \mathbb{R}^d \to \mathbb{R}^p$, and if there exists a gauge-integrable function $h : \mathbb{R}^d \to \mathbb{R}$ such that $|f_n| \leq h$ in $\mathbb{R}^d$ for every $n \in \mathbb{N}$, then $f$ is gauge-integrable and*

$$\lim_{n \to \infty} \int_{\mathbb{R}^d} f_n = \int_{\mathbb{R}^d} f.$$

*Proof of Theorem 1.3.* If $f$ is gauge-integrable, then $f$ is gauge-measurable by Proposition 5.2, and it follows from the Absolute Cauchy criterion above that the function $h := |f|$ is gauge-integrable.

Conversely, if $f$ is gauge-measurable, then by Proposition 2.2, for every $n \in \mathbb{N}$ the truncated function $T_n \circ f$ is also gauge-measurable, where $T_n : \mathbb{R}^p \to \mathbb{R}^p$ is the truncation function defined for $w \in \mathbb{R}^p$ by

$$T_n(w) = \begin{cases} w & \text{if } |w| \leq n, \\ nw/|w| & \text{if } |w| > n, \end{cases}$$

Since the function $T_n \circ f$ is bounded and satisfies $|T_n \circ f| \leq |f| \leq h$ in $\mathbb{R}^d$, it follows from Proposition 5.3 that $T_n \circ f$ is gauge-integrable, and we conclude applying the Dominated convergence theorem for gauge integrals as $n$ tends to infinity. □

## 6 Pointwise limit

A crucial feature of Lebesgue-measurable functions is its stability under pointwise convergence. Up to now, we only have proved the stability of gauge measurability under uniform convergence (Proposition 2.3). Thanks to the relationship that we have established between gauge measurability and gauge integrability, we now obtain a pointwise convergence property.

**Proposition 6.1** (Pointwise limit). *Let $(f_n)_{n \in \mathbb{N}}$ be a sequence of gauge-measurable functions from $\mathbb{R}^d$ to $\mathbb{R}^p$. If $(f_n)_{n \in \mathbb{N}}$ converges pointwise to the function $f : \mathbb{R}^d \to \mathbb{R}^p$, then $f$ is gauge-measurable.*

We first prove two particular cases of this proposition, which as we shall see yield the general case. We denote the *characteristic function* of a set $A \subset \mathbb{R}^d$ by $\chi_A$, that is $\chi_A : \mathbb{R}^d \to \mathbb{R}$ is the function defined for each $x \in \mathbb{R}^d$ by

$$\chi_A(x) = \begin{cases} 1 & \text{if } x \in A, \\ 0 & \text{if } x \notin A. \end{cases}$$

**Lemma 6.2.** *Let $(A_l)_{l \in \mathbb{N}}$ be an increasing sequence of open subsets which cover $\mathbb{R}^d$ with $\overline{A_{l-1}} \subset A_l$ for every $l \in \mathbb{N}_*$. If a function $f : \mathbb{R}^d \to \mathbb{R}^p$ is such that $f\chi_{A_l}$ is gauge-measurable for every $l \in \mathbb{N}$, then $f$ is also gauge-measurable.*

*Proof.* Given $l, m \in \mathbb{N}_*$ with $l \leq m$, we first observe that if

$$(A_l \setminus A_{l-2}) \cap (A_m \setminus A_{m-2}) \neq \emptyset, \tag{6.1}$$

then by monotonicity of the sequence $(A_l)_{l \in \mathbb{N}}$ we have $m = l$ or $m = l + 1$. Here, we use the convention that $A_{-1} = \emptyset$. Now let $(\gamma_l)_{l \in \mathbb{N} \setminus \{0,1\}}$ be a sequence of gauges on $\mathbb{R}^d$ to be chosen later on. We define a new gauge $\gamma$ on $\mathbb{R}^d$ as follows: for every $x \in \mathbb{R}^d$, denote by $l$ the smallest integer in $\mathbb{N}_*$ such that $x \in A_l$ and let

$$\gamma(x) = \gamma_l(x) \cap \gamma_{l+1}(x) \cap (A_l \setminus \overline{A_{l-2}}).$$

Since $\overline{A_{l-2}} \subset A_{l-1}$ by assumption, we have $x \notin \overline{A_{l-2}}$, and then the open set $A_l \setminus \overline{A_{l-2}}$ contains $x$. Thus, $\gamma$ is a well-defined gauge on $\mathbb{R}^d$.

Given $\eta > 0$, we claim that

$$U_{\gamma, \eta} \subset \bigcup_{l \in \mathbb{N}_*} V_{l+1}, \tag{6.2}$$

where

$$V_l := \bigcup_{\substack{x,z \in \mathbb{R}^d \\ |f\chi_{A_l}(x) - f\chi_{A_l}(z)| \geq \eta}} (\gamma_l(x) \cap \gamma_l(z)).$$

Indeed, assume that $x, z \in \mathbb{R}^d$ are such that $|f(x) - f(z)| \geq \eta$. Let $l$ and $m$ be the smallest integers in $\mathbb{N}_*$ such that $x \in A_l$ and $z \in A_m$; we may assume without loss of generality that $l \leq m$. If $\gamma(x) \cap \gamma(z) \neq \emptyset$, then (6.1) holds, and thus $m = l$ or $m = l + 1$. Hence, we have

$$|f\chi_{A_{l+1}}(x) - f\chi_{A_{l+1}}(z)| = |f(x) - f(z)| \geq \eta$$

and

$$\gamma(x) \cap \gamma(z) \subset \gamma_{l+1}(x) \cap \gamma_{l+1}(z) \subset V_{l+1},$$

which implies (6.2).

Let $\epsilon > 0$. Since the function $f\chi_{A_l}$ is gauge-measurable, by Lemma 4.2 we can choose the gauge $\gamma_l$ on $\mathbb{R}^d$ such that $\mu(V_l) \leq \epsilon/2^{l-1}$. Thus, by the inclusion (6.2) and the countable subadditivity of the inner measure $\mu$ we get

$$\mu(U_{\gamma, \eta}) \leq \sum_{l \in \mathbb{N}_*} \mu(V_{l+1}) \leq \sum_{l \in \mathbb{N}_*} \frac{\epsilon}{2^l} = \epsilon.$$

By Lemma 4.2, we deduce that $f$ is gauge-measurable. □

**Lemma 6.3.** *If the function $f : \mathbb{R}^d \to \mathbb{R}^p$ is such that the truncation $T_j \circ f$ is gauge-measurable for every $j \in \mathbb{N}$, then $f$ is gauge-measurable.*

*Proof.* Given a sequence of gauges $(\gamma_j)_{j \in \mathbb{N}_*}$ on $\mathbb{R}^d$, consider the gauge $\gamma$ defined for $x \in \mathbb{R}^d$ by

$$\gamma(x) = \gamma_0(x) \cap \cdots \cap \gamma_{j+1}(x),$$

where $j \in \mathbb{N}$ is the smallest integer such that $|f(x)| \leq j$. For every $0 < \eta \leq 1$, we claim that
$$U_{\gamma,\eta} \subset \bigcup_{j \in \mathbb{N}} W_{j+1}, \tag{6.3}$$
where
$$W_j := \bigcup_{\substack{x,z \in \mathbb{R}^d \\ |T_j \circ f(x) - T_j \circ f(z)| \geq \eta}} \bigl(\gamma_j(x) \cap \gamma_j(z)\bigr).$$

For this purpose, for every $x, z \in \mathbb{R}^d$ such that $|f(x) - f(z)| \geq \eta$, which we may assume that $|f(z)| \geq |f(x)|$, let $j \in \mathbb{N}$ be the smallest integer such that $|f(x)| \leq j$. Since $\eta \leq 1$, we also have
$$|T_{j+1} \circ f(x) - T_{j+1} \circ f(z)| \geq \eta.$$
From the choice of the gauge $\gamma$, we deduce that
$$\gamma(x) \cap \gamma(z) \subset \gamma_{j+1}(x) \cap \gamma_{j+1}(z) \subset W_{j+1},$$
and the inclusion (6.3) follows.

Let $\epsilon > 0$. Since the function $T_j \circ f$ is gauge-measurable, by Lemma 4.2 we can choose the gauge $\gamma_j$ on $\mathbb{R}^d$ such that $\mu(W_j) \leq \epsilon/2^j$. Proceeding as in the previous lemma, we have $\mu(U_{\gamma,\eta}) \leq \epsilon$, hence $f$ is gauge-measurable. □

*Proof of Proposition 6.1.* We first assume that there exists a gauge-integrable function $h : \mathbb{R}^d \to \mathbb{R}$ such that $|f_n| \leq h$ in $\mathbb{R}^d$ for every $n \in \mathbb{N}$. By Theorem 1.3, each function $f_n$ is gauge-integrable, and it then follows from the Dominated convergence theorem that $f$ is gauge-integrable, hence also gauge-measurable.

In the general case where the sequence $(f_n)_{n \in \mathbb{N}}$ need not be bounded by an integrable function, for every $n, l, j \in \mathbb{N}$ we consider the function
$$g_{n,l,j} = (T_j \circ f_n)\chi_{B_{l+1}(0)}.$$

These functions are all gauge-measurable. Indeed, $T_j \circ f_n$ is gauge-measurable by composition with the uniformly continuous function $T_j$ (Proposition 2.2), and thus $g_{n,l,j}$ is gauge-measurable as the product of bounded gauge-measurable functions (Proposition 3.2).

Since $|g_{n,l,j}| \leq j\chi_{B_{l+1}(0)}$ in $\mathbb{R}^d$ and the characteristic function $\chi_{B_{l+1}(0)}$ is gauge-integrable, as $n$ tends to infinity it follows from the first case we considered above that the functions $(T_j \circ f)\chi_{B_{l+1}(0)}$ are gauge-measurable for every $l, j \in \mathbb{N}$. By Lemma 6.2, as $l$ tends to infinity we deduce that $T_j \circ f$ is gauge-measurable for every $j \in \mathbb{N}$. The conclusion then follows from Lemma 6.3 as $j$ tends to infinity. □

A consequence of Propositions 3.2 and 6.1 is that the product of gauge-measurable functions is also gauge-measurable:

**Proposition 6.4** (Product). *If the functions $f : \mathbb{R}^d \to \mathbb{R}^p$ and $g : \mathbb{R}^d \to \mathbb{R}$ are gauge-measurable, then their product $fg$ is also gauge-measurable.*

More generally, we can weaken the assumptions of Proposition 2.2 on the gauge measurability of composite functions:

**Proposition 6.5** (Composition with continuous functions). *If the function $f : \mathbb{R}^d \to \mathbb{R}^p$ is gauge-measurable and the function $\Phi : \mathbb{R}^p \to \mathbb{R}^\ell$ is continuous, then the composition $\Phi \circ f : \mathbb{R}^d \to \mathbb{R}^\ell$ is gauge-measurable.*

*Proof.* We consider a continuous function $\varphi : \mathbb{R}^p \to \mathbb{R}$ with compact support, and we define $\Phi_n : \mathbb{R}^p \to \mathbb{R}^\ell$ for each $n \in \mathbb{N}_*$ and $y \in \mathbb{R}^p$ by $\Phi_n(y) = \varphi(y/n)\Phi(y)$. Since the function $\Phi_n$ is continuous and has compact support, $\Phi_n$ is uniformly continuous, and thus, in view of Proposition 2.2, the function $\Phi_n \circ f$ is gauge-measurable. We conclude by observing that, for every $x \in \mathbb{R}^n$, the sequence $(\Phi_n(f(x)))_{n \in \mathbb{N}_*}$ converges to $\Phi(f(x))$ provided that $\varphi(0) = 1$, and thus by Proposition 6.1 the function $\Phi \circ f$ is gauge-measurable. □

The proof of Proposition 6.5 shows that the class of functions $\Phi : \mathbb{R}^p \to \mathbb{R}^\ell$ such that for every gauge-measurable function $f : \mathbb{R}^d \to \mathbb{R}^p$, the composition $\Phi \circ f$ is measurable is stable under pointwise convergence. This class thus forms a Baire system and contains in particular all Baire (or analytic representable) functions, which coincide by the Lebesgue–Hausdorff theorem with all Borel-measurable functions [9, Theorem 43.IV; 12, §31].

Another consequence of Proposition 6.1 combined with the gauge measurability of gauge-integrable functions (Proposition 5.2) is that the pointwise limit of a sequence of gauge-integrable functions is always gauge-measurable. This implies in particular that measurable functions in the sense of McShane [23, Definition III-10-1] are indeed gauge-measurable. Conversely, every gauge-measurable function $f : \mathbb{R}^d \to \mathbb{R}^p$ in the sense of Definition 1.2 is the limit of a sequence of gauge-integrable functions. This assertion follows from a diagonalization procedure using the functions $g_{n,l,j}$ which are used in the proof of Proposition 6.1 above. For example, the sequence of gauge-integrable functions $(g_{n,n,n})_{n \in \mathbb{N}}$ converges pointwise to the measurable function $f$. Another pointwise approximation of $f$ in terms of gauge-measurable step functions is pursued in Section 8.

## 7 Gauge-measurable sets

We define gauge measurability of a set in the spirit of its counterpart for functions:

**Definition 7.1.** *A set $A \subset \mathbb{R}^d$ is gauge-measurable whenever, for every $\epsilon > 0$, there exists a gauge $\gamma$ on $\mathbb{R}^d$ such that, for every finite set of disjoint rectangles $(R_i)_{i \in \{1,\dots,k\}}$, every finite set of points $(c_i)_{i \in \{1,\dots,k\}}$ contained in $A$, and every finite set of points $(c'_i)_{i \in \{1,\dots,k\}}$ contained in $\mathbb{R}^d \setminus A$ that satisfy $R_i \subset \gamma(c_i) \cap \gamma(c'_i)$ for every $i \in \{1,\dots,k\}$, one has*

$$\sum_{i=1}^{k} \mathrm{vol}\,(R_i) \le \epsilon.$$

It follows from this definition that $A$ is gauge-measurable if and only if its complement $\mathbb{R}^d \setminus A$ is gauge-measurable. Also observe that for any $0 < \eta < 1$ we have
$$|\chi_A(x) - \chi_A(z)| \geq \eta$$
if and only if $x \in A$ and $z \in \mathbb{R}^d \setminus A$, or $x \in \mathbb{R}^d \setminus A$ and $z \in A$. In view of Definitions 1.2 and 7.1, it thus follows that the set $A \subset \mathbb{R}^d$ is gauge-measurable if and only if the characteristic function $\chi_A$ is gauge-measurable.

As in Lemma 4.2, the definition above can be reformulated by replacing the tagged partitions with the inner measure of an open set:

**Lemma 7.2** (Gauge-intersection characterization). *The set $A \subset \mathbb{R}^d$ is gauge-measurable if and only if, for every $\epsilon > 0$, there exists a gauge $\gamma$ on $\mathbb{R}^d$ such that the open set*
$$U_{A,\gamma} := \bigcup_{\substack{x \in A, \\ z \in \mathbb{R}^d \setminus A}} \bigl(\gamma(x) \cap \gamma(z)\bigr)$$
*satisfies $\mu(U_{A,\gamma}) \leq \epsilon$.*

This characterization can be established along the lines of the proof of Lemma 4.2. The family of gauge-measurable sets forms an algebra:

**Proposition 7.3.** *If the sets $A_1, A_2 \subset \mathbb{R}^d$ are gauge-measurable, then $A_1 \cup A_2$, $A_1 \cap A_2$, and $A_1 \setminus A_2$ are also gauge-measurable.*

*Proof.* We prove that $A_1 \cup A_2$ is gauge-measurable. For this purpose, observe that every $z \in \mathbb{R}^d \setminus (A_1 \cup A_2)$ satisfies $z \in \mathbb{R}^d \setminus A_1$ and $z \in \mathbb{R}^d \setminus A_2$. Thus, given a gauge $\gamma$ on $\mathbb{R}^d$, we have
$$U_{A_1 \cup A_2, \gamma} = U_{A_1, \gamma} \cup U_{A_2, \gamma}$$

Given $\epsilon > 0$, let $\gamma_1$ and $\gamma_2$ be two gauges on $\mathbb{R}^d$ satisfying the conclusion of Lemma 7.2 for $A_1$ and $A_2$, respectively, with parameter $\epsilon/2$. Take the gauge $\gamma$ defined for $x \in \mathbb{R}^d$ by $\gamma(x) = \gamma_1(x) \cap \gamma_2(x)$. Thus,
$$U_{A_1 \cup A_2, \gamma} = U_{A_1, \gamma} \cup U_{A_2, \gamma} \subset U_{A_1, \gamma_1} \cup U_{A_2, \gamma_2},$$
and by subadditivity of $\mu$ we then get
$$\mu(U_{A_1 \cup A_2, \gamma}) \leq \mu(U_{A_1, \gamma_1}) + \mu(U_{A_2, \gamma_2}) \leq \frac{\epsilon}{2} + \frac{\epsilon}{2} = \epsilon.$$
Hence, $A_1 \cup A_2$ is gauge-measurable.

Since we have
$$\mathbb{R}^d \setminus (A_1 \cap A_2) = (\mathbb{R}^d \setminus A_1) \cup (\mathbb{R}^d \setminus A_2)$$
and both sets $\mathbb{R}^d \setminus A_1$ and $\mathbb{R}^d \setminus A_2$ are gauge-measurable, we deduce that $\mathbb{R}^d \setminus (A_1 \cap A_2)$ is gauge-measurable, and thus the intersection $A_1 \cap A_2$ is also gauge-measurable. Finally, since
$$A_1 \setminus A_2 = A_1 \cap (\mathbb{R}^d \setminus A_2)$$
is the intersection of two gauge-measurable sets, $A_1 \setminus A_2$ is also gauge-measurable. $\square$

Using the equivalence between the gauge measurability of the set $A$ and the gauge-measurability of the characteristic function $\chi_A$, we deduce that the family of gauge-measurable sets forms a $\sigma$-algebra:

**Proposition 7.4** (Countable union). *If $(A_n)_{n\in\mathbb{N}}$ is a sequence of gauge-measurable sets in $\mathbb{R}^d$, then the set $\bigcup_{k\in\mathbb{N}} A_k$ is also gauge-measurable.*

*Proof.* The sequence of characteristic functions $(f_n)_{n\in\mathbb{N}}$ defined for each $n \in \mathbb{N}$ by $f_n = \chi_{\bigcup_{k=0}^n A_k}$ converges pointwise to the characteristic function $\chi_{\bigcup_{k\in\mathbb{N}} A_k}$ in $\mathbb{R}^d$. By induction using Proposition 7.3, each set $\bigcup_{k=0}^n A_k$ is gauge-measurable and thus each function $f_n$ is gauge-measurable. From the stability of gauge measurability under pointwise convergence (Proposition 6.1), we deduce that the function $\chi_{\bigcup_{k\in\mathbb{N}} A_k}$ is also gauge-measurable, hence the set $\bigcup_{k\in\mathbb{N}} A_k$ is gauge-measurable. □

Let us now prove Lebesgue's regularity property which provides one with a necessary and sufficient condition for a set to be gauge-measurable, and also implies that a set is gauge-measurable if and if it is Lebesgue-measurable [26, Lemma 3.22].

**Proposition 7.5** (Regularity). *The set $A \subset \mathbb{R}^d$ is gauge-measurable if and only if, for every $\epsilon > 0$, there exist an open set $V \subset \mathbb{R}^d$ and a closed set $C \subset \mathbb{R}^d$ such that $C \subset A \subset V$ and $\mu(V \setminus C) \leq \epsilon$.*

*Proof.* Given a gauge $\gamma$ on $\mathbb{R}^d$, set

$$V = \bigcup_{x \in A} \gamma(x) \quad \text{and} \quad C = \bigcap_{z \in \mathbb{R}^d \setminus A} \bigl(\mathbb{R}^d \setminus \gamma(z)\bigr).$$

Observe that $V$ is open, $C$ is closed, and $V \setminus C = U_{A,\gamma}$. Thus, given $\epsilon > 0$, if the set $A$ is gauge-measurable and one takes a gauge $\gamma$ such that $\mu(U_{A,\gamma}) \leq \epsilon$, then the sets $V$ and $C$ above satisfy the requirements.

Conversely, if the set $A$ satisfies the regularity condition, then given $\epsilon > 0$ we take the sets $V$ and $C$ as in the statement. The gauge $\gamma$ defined on $\mathbb{R}^d$ by setting $\gamma(x) = V$ if $x \in A$ and $\gamma(x) = \mathbb{R}^d \setminus C$ if $x \in \mathbb{R}^d \setminus A$ satisfies $U_{A,\gamma} = V \setminus C$, and thus $\mu(U_{A,\gamma}) \leq \epsilon$. Hence, the set $A$ is gauge-measurable by Lemma 7.2. □

## 8 Pointwise approximation

We conclude with the pointwise approximation of a gauge-measurable function by step functions. We recall that $g : \mathbb{R}^d \to \mathbb{R}^p$ is a step function if the image $g(\mathbb{R}^d)$ is a finite set.

**Proposition 8.1.** *If a function $f : \mathbb{R}^d \to \mathbb{R}^p$ is gauge-measurable, then there exists a sequence of gauge-measurable step functions $(f_n)_{n\in\mathbb{N}}$ from $\mathbb{R}^d$ to $\mathbb{R}^p$ which converges pointwise to $f$ in $\mathbb{R}^d$ and satisfies $|f_n| \leq |f|$ in $\mathbb{R}^d$ for every $n \in \mathbb{N}$.*

This statement allows one to recover the usual strategy to define the Lebesgue integral via measurable step functions. In our case, if $f$ is gauge-integrable, and thus $|f|$ is also gauge-integrable by the Absolute Cauchy criterion, then from the Dominated convergence theorem (Proposition 5.4) we indeed have that

$$\int_{\mathbb{R}^d} f = \lim_{n \to \infty} \int_{\mathbb{R}^d} f_n.$$

The difference here is that this is a *property* of the gauge integral, rather than a definition.

Before proving Proposition 8.1, we first study the inverse image of rectangles by gauge-measurable functions:

**Proposition 8.2.** *If a function $f : \mathbb{R}^d \to \mathbb{R}^p$ is gauge-measurable, then, for every rectangle $R \subset \mathbb{R}^p$, the set $f^{-1}(R)$ is gauge-measurable.*

*Proof.* Observe that $\chi_R \circ f = \chi_{f^{-1}(R)}$. To prove the proposition, it thus suffices to prove that the function $\chi_R \circ f$ is gauge-measurable. For this purpose, take a sequence of uniformly continuous functions $(\Phi_n)_{n \in \mathbb{N}}$ from $\mathbb{R}^p$ to $\mathbb{R}$ which converges pointwise to $\chi_R$. Then, by Proposition 2.2 the function $\Phi_n \circ f$ is gauge-measurable for every $n \in \mathbb{N}$, and the sequence $(\Phi_n \circ f)_{n \in \mathbb{N}}$ converges pointwise to $\chi_R \circ f$. By the stability property of sequences of gauge-measurable functions (Proposition 6.1), we deduce that $\chi_R \circ f$ is gauge-measurable, and the conclusion follows. □

The converse of Proposition 8.2 is also true: if $f^{-1}(R)$ is gauge-measurable for every rectangle $R \subset \mathbb{R}^p$, then $f$ is gauge-measurable. This assertion can be deduced from the proof of Proposition 8.1 below, since under such an assumption the functions $f_n$ which are defined in (8.1) below are all gauge-measurable and the function $f$ is the pointwise limit of the sequence $(f_n)_{n \in \mathbb{N}}$.

*Proof of Proposition 8.1.* Take a sequence of positive numbers $(\epsilon_n)_{n \in \mathbb{N}}$ that converges to 0. For each $n \in \mathbb{N}$, let $(R_{i,n})_{i \in \{1,\dots,k_n\}}$ be a finite family of disjoint rectangles whose diameters do not exceed $\epsilon_n$ that covers the ball $B_{n+1}(0)$ in $\mathbb{R}^p$. For each $i \in \{1, \dots, k_n\}$, let $a_i$ be a point with smallest norm in $\overline{R_{i,n}}$, and define

$$f_n := \sum_{i=1}^{k_n} a_i \, \chi_{f^{-1}(R_{i,n})}. \tag{8.1}$$

For every $x \in f^{-1}(B_{n+1}(0))$, we then have

$$|f_n(x) - f(x)| \leq \epsilon_n,$$

hence the sequence $(f_n)_{n \in \mathbb{N}}$ converges pointwise to $f$ in $\mathbb{R}^d$. [The convergence is uniform when $f$ is a bounded function.] By the choice of the point $a_i$, we also have

$$|f_n(x)| \leq |f(x)|$$

if $f(x) \in \bigcup_{i=1}^{k_n} R_{i,n}$, while the left-hand side vanishes otherwise. This estimate thus holds for every $x \in \mathbb{R}^d$.

Assuming that the set $f^{-1}(R)$ is gauge-measurable for every rectangle $R \subset \mathbb{R}^p$, which by Proposition 8.2 is the case when the function $f$ is gauge-measurable, it follows from the linear stability of gauge-measurable functions (Proposition 3.1) that $f_n$ is a gauge-measurable step function, and this gives the conclusion.   □

# References

[1] R. G. Bartle, *A modern theory of integration*, Graduate Studies in Mathematics, vol. 32, American Mathematical Society, Providence, RI, 2001. ↑2

[2] J.-M. Bony, G. Choquet, and G. Lebeau, *Le centenaire de l'intégrale de Lebesgue*, C. R. Acad. Sci. Paris Sér. I Math. **332** (2001), 85–90. ↑1

[3] N. Bourbaki, *Éléments de mathématique. Fasc. XIII. Livre VI: Intégration. Chapitres 1, 2, 3 et 4: Inégalités de convexité, Espaces de Riesz, Mesures sur les espaces localement compacts, Prolongement d'une mesure, Espaces $L^p$*, 2nd ed., Actualités Scientifiques et Industrielles, No. 1175, Hermann, Paris, 1965. ↑7

[4] A.-L. Cauchy, *Résumé des leçons données à l'École royale polytechnique sur le calcul infinitésimal*, Vol. I, Imprimerie royale, Debure frères, Paris, 1823. ↑1

[5] Ch.-J. de la Vallée Poussin, *Cours d'analyse infinitésimale*, 2nd ed., Vol. I, Uyspruyt, Louvain, 1909. ↑2

[6] \_\_\_\_\_\_, *Cours d'analyse infinitésimale*, 2nd ed., Vol. II, Uyspruyt, Louvain, 1912. ↑2

[7] \_\_\_\_\_\_, *Sur l'intégrale de Lebesgue*, Trans. Amer. Math. Soc. **16** (1915), 435–501. ↑2

[8] A. Denjoy, *Une extension de l'integrale de M. Lebesgue*, C. R. Acad. Sci. Paris **154** (1912), 859–862. ↑2

[9] F. Hausdorff, *Set theory*, 2nd ed., Chelsea Publishing, New York, 1962. ↑16

[10] R. Henstock, *Definitions of Riemann type of the variational integrals*, Proc. London Math. Soc. (3) **11** (1961), 402–418. ↑1

[11] J.-P. Kahane, *Naissance et postérité de l'intégrale de Lebesgue*, Gaz. Math. **89** (2001), 5–20. ↑1

[12] K. Kuratowski, *Topology*, Vol. I, Academic Press, New York-London, 1966. ↑16

[13] D. S. Kurtz and C. W. Swartz, *Theories of integration. The integrals of Riemann, Lebesgue, Henstock-Kurzweil, and McShane*, 2nd ed., Series in Real Analysis, vol. 13, World Scientific Publishing Co. Pte. Ltd., Hackensack, NJ, 2012. ↑10

[14] J. Kurzweil, *Generalized ordinary differential equations and continuous dependence on a parameter*, Czechoslovak Math. J. **7 (82)** (1957), 418–449 (Russian). ↑1

[15] H. Lebesgue, *Sur une généralisation de l'intégrale définie*, C. R. Acad. Sci. Paris **132** (1901), 1025–1028. ↑1

[16] T. Y. Lee, *Henstock-Kurzweil integration on Euclidean spaces*, Series in Real Analysis, vol. 12, World Scientific Publishing, Hackensack, NJ, 2011. ↑2

[17] N. Lusin, *Sur les propriétés de l'intégrale de M. Denjoy*, C. R. Acad. Sci. Paris **155** (1912), 1475–1478. ↑2

[18] J. Mawhin, *Introduction à l'analyse*, Cabay, Louvain, 1979. ↑2, 3

[19] \_\_\_\_\_\_, *Analyse: Fondements, techniques, évolution*, 2nd ed., Accès Sciences, De Boeck, Bruxelles, 1997. ↑2, 3

[20] \_\_\_\_\_\_, *Two histories of integration theory: Riemannesque vs Romanesque*, Acad. Roy. Belg. Bull. Cl. Sci. (6) **18** (2007), 47–63. ↑1

[21] \_\_\_\_\_\_, *The* Cours d'analyse infinitésimale *of Charles-Jean de La Vallée Poussin: from innovation to tradition*, Jahresber. Dtsch. Math.-Ver. **116** (2014), 243–259. ↑2

[22] E. J. McShane, *A unified theory of integration*, Amer. Math. Monthly **80** (1973), 349–359. ↑1

[23] \_\_\_\_\_\_, *Unified integration*, Pure and Applied Mathematics, vol. 107, Academic Press, New York, 1983. ↑1, 2, 4, 10, 12, 16

[24] O. Perron, *Über den Integralbegriff*, S.-B. Heidelberg Akad. Wiss., Abt. A **16** (1914), 1–16. ↑2

[25] B. Riemann, *Über die Darstellbarkeit einer Function durch eine trigonometrische Reihe*, Abhandlungen der Königlichen Gesellschaft der Wissenschaften in Göttingen **13** (1868), 87–131. ↑1

[26] R. L. Wheeden and A. Zygmund, *Measure and integral. An introduction to Real Analysis*, 2nd ed., Pure and Applied Mathematics, CRC Press, Boca Raton, FL, 2015. ↑18